\documentclass[a4paper,10pt]{article}
\usepackage[utf8]{inputenc}
\usepackage[english]{babel}
\usepackage{amssymb, amsmath, enumerate}
\usepackage[all,2cell]{xy}
\usepackage{hyperref}
\usepackage{bm}

\usepackage{theorem}
\newtheorem{teo}{Theorem}[section]
\newtheorem{lem}[teo]{Lemma} 

\newtheorem{prop}[teo]{Proposition} 

\theorembodyfont{\normalfont\upshape} 
\newtheorem{defn}[teo]{Definition} 
\newtheorem{ex}[teo]{Example}

\newenvironment{dem}[1][Proof]{\textbf{#1:}\

}  {\hfill\rule{1ex}{1ex}}

\newcommand{\PQM}{{\bf PQM}}
\newcommand{\PQMrr}{{\bf PQMrr}}
\newcommand{\pPQM}{{\bf pPQM}}
\newcommand{\pPQMrr}{{\bf pPQMrr}}

\author{
  Hugo Rafael Ribeiro\\
  \texttt{hugor@ime.usp.br}
  \and
  Kaique Matias de Andrade Roberto\\
  \texttt{kaique.roberto@usp.br}
  \and
  Hugo Luiz Mariano\\
  \texttt{hugomar@ime.usp.br} \\
  Instituto de Matematica e Estatistica \\
  Universidade de Sao Paulo, Brazil
}
\title{Promiscuously Quadratic Rings}

\begin{document}

\maketitle

\begin{abstract}
We consider a convenient category of ``quadratic'' multirings, that allows simple functorial relations 
with categories associated with abstract quadratic forms theories and shares many good aspects of the 
theories of Special Groups and of Realsemigroups, but tries to avoids its drawbacks.
\end{abstract}

\section{Introduction}
\hspace*{\parindent}
The main motivation for this work is to determine a simple (sub)category of multirings that:

\begin{enumerate}[I -]
 \item shares many good aspects of the theories of Special Groups and of RealSemigroups, but tries to 
avoids its drawbacks;

 \item have simple functorial relations with the categories of Special Groups and of Realsemigroups;

 \item allows a natural approach to the representation problem in the abstract quadratic forms theories 
of Special Groups and of Realsemigroups.
\end{enumerate}

Concerning Special Groups and of Realsemigroups:

\begin{itemize}
 \item Both shares the nice property  of being first-order theories, axiomatized by simple axioms.
 \item The theory of Special Groups allows a nice n-dimensional isometry theory, intrinsically related 
with a ring. On the other hand , being a  quadratic form theory obtained only by the invertible 
coefficient of a ring, it  loses on the spot many natural rings as $\mathbb{Z}$.
 \item The theory of Realsemigroups presents an approach of quadratic form theory that allows consider 
general coefficients in a ring. In this setting, on other hand,  there is no  intrinsic definition of 
isometry or a nice n-dimensional isometry theory.
\end{itemize}

\section{Basic Concepts}
\hspace*{\parindent}
This section is a compiled o basic definitions and results included for the convenience of 
the reader. For more details and examples, consult \cite{marshall2006real}, \cite{dickmann2000special}, 
\cite{dickmann2004real} or \cite{ribeiro2016functorial}.

\begin{defn}\label{defn:multigroupI}
 A multigroup is a quadruple $(G,\ast,r,1)$, where $G$ is a non-empty set, \newline $\ast:G\times 
G\rightarrow\mathcal P(G)\setminus\{\emptyset\}$ and
 $r:G\rightarrow G$
 are functions, and $1$ is an element of $G$ satisfying:
 \begin{enumerate}[i -]
  \item If $ z\in x\ast y$ then $x\in z\ast r(y)$ and $y\in r(x)\ast z$.
  \item $y\in 1\ast x$ iff $x=y$.
  \item With the convention $x\ast(y\ast z)=\bigcup\limits_{w\in y\ast z}x\ast w$ and 
  $(x\ast y)\ast z=\bigcup\limits_{t\in x\ast y}t\ast z$,
  $$x\ast(y\ast z)=(x\ast y)\ast z\mbox{for all }x,y,z\in G.$$
  
	A multigroup is said to be \textit{commutative} if
  \item $x\ast y=y\ast x$ for all $x,y\in G$.
 \end{enumerate}
\end{defn}

\begin{defn}\label{defn:multiring}
 A multiring is a sextuple $(R,+,\cdot,-,0,1)$ where $R$ is a non-empty set, $+:R\times 
R\rightarrow\mathcal P(R)\setminus\{\emptyset\}$,
 $\cdot:R\times R\rightarrow R$
 and $-:R\rightarrow R$ are functions, $0$ and $1$ are elements of $R$ satisfying:
 \begin{enumerate}[i -]
  \item $(R,+,-,0)$ is a commutative multigroup;
  \item $(R,\cdot,1)$ is a commutative monoid;
  \item $a0=0$ for all $a\in R$;
  \item If $c\in a+b$, then $cd\in ad+bd$. Or equivalently, $(a+b)d\subseteq ab+bd$.
 \end{enumerate}
 
 $R$ is said to be a multidomain if do not have zero divisors, and $R$ will be a multifield if every 
non-zero element of $R$ has multiplicative inverse. If $(a+b)d=ab+bd$ for all $a,b,d\in R$, then $R$ 
will be an {\em hyperring}.
\end{defn}

\begin{ex}\label{ex:1.3}
$ $
 \begin{enumerate}[a -]
  \item Every ring, domain and field is a multiring, multidomain and multifield respectively.
  
  \item $Q_2=\{-1,0,1\}$ is a multifield with the usual product and the multivalued sum defined by 
relations
  $$\begin{cases}
     0+x=x+0=x,\,\mbox{for every }x\in Q_2 \\
     1+1=1,\,(-1)+(-1)=-1 \\
     1+(-1)=(-1)+1=\{-1,0,1\}
    \end{cases}
  $$

  \item In the set $\mathbb R_+$ of positive real numbers, we define $a\bigtriangledown 
b=\{c\in\mathbb R_+:|a-b|\le 
c\le a+b\}$. 
We have that
  $\mathbb R_+$ with the usual product and $\bigtriangledown$ multivalued sum is a multifield, called 
triangle 
multifield 
\cite{viro2010hyperfields}. We
  denote this multifield by $\mathcal{T}\mathbb R_+$.
\end{enumerate}
\end{ex}

\begin{lem}\label{lemma:1.4}
 Let $F$ be a multifield. Then $(a+b)d=ad+bd$ for every $a,b,d\in F$.
\end{lem}
\begin{dem}
 We have $(a+b)d\subseteq ad+bd$ already. For the other inclusion, if $d=0$, it is done. If $d\ne0$, we 
have:
 \begin{align*}(ad+bd)d^{-1}\subseteq (ad)d^{-1}+(bd)d^{-1}=ad+bd\Rightarrow \\ 
  ad+bd=[(ad+bd)d^{-1}]d\subseteq(a+b)d.\end{align*}
\end{dem}
   
  Ideals, spectrum, orderings and another constructions in the multivalued language can be found in 
\cite{marshall2006real}. 
   
\begin{defn}\label{defn:morphism}
 Let $A$ and $B$ multirings. A map $f:A\rightarrow B$ is a morphism if for all $a,b,c\in A$:
 \begin{enumerate}[i -]
  \item $c\in a+b\Rightarrow f(c)\in f(a)+f(b)$;
  \item $f(-a)=-f(a)$;
  \item $f(0)=0$;
  \item $f(ab)=f(a)f(b)$;
  \item $f(1)=1$.
 \end{enumerate}
\end{defn}

 For multirings, there are various sorts of ``substrucutre'' that one can consider. If $A,B$ 
are multirings, we say $A$ is embedded in $B$ by the morphism $\iota:A\rightarrow B$ if $\iota$ 
is injective. We say $A$ is strongly embedded in $B$ if $A$ is
 embedded in $B$ and, for all $a,b,c\in A$, $\iota(c)\in\iota(a)+_B\iota(b)\Rightarrow c\in a+_Ab$. We 
say $A$ is a submultiring of $B$ if $A$
 is strongly embedded in $B$ and, for all $a,b\in A$ and all $c\in B$, 
$c\in\iota(a)+_B\iota(b)\Rightarrow c\in\iota(A)$. Note that in the rings case, these 
all definitions coincide.
 
 The category of multifields (respectively multirings) and theirs morphisms will be 
denoted by $\mathcal{MF}$ (respectively $\mathcal{MR}$). Now, we present the main construction:
 
 \begin{defn}[Marshall's Quotient]\label{defn:strangeloc}
 Fix a multiring $A$ and a multiplicative subset $S$ of $A$. Define an equivalence relation $\sim$ 
on $A$ by $a\sim b$ iff $as=bt$ for some $s,t\in S$. Denote by $\overline a$ the equivalence class of 
$a$ and set $A/_mS=\{\overline a:a\in A\}$. Defining $\overline a+\overline b=\{\overline c:cv\in 
as+bt,\,\mbox{for some }s,t,v\in S\}$, $-\overline 
a=\overline{-a}$, and $\overline{a}\overline{b}=\overline{ab}$ we have that 
$(A/_mS,+,\cdot,-,\overline{0},\overline{1})$ is a multiring. When $S=\sum A^{*2}$, we will denote 
$A/_m\sum A^{*2}=Q_{\mbox{red}(A)}$.
 \end{defn}

 \begin{defn}[Corollary 4.2 in \cite{marshall2006real}]\label{defn:mfrealreduced}
 A multifield $F$ is said to be \textit{real reduced} if $a^3=a$ and $(a\in1+1)\Rightarrow(a=1)$ for all 
$a\in F$. This implies that the morphism $a\mapsto\overline a$ from $F$ to $Q_{\mbox{red}}(F)$ is an 
isomorphism.

A morphism of real reduced multifield is just a morphism of multifields. The category of real 
reduced multifields will be denoted by $\mathcal{MF}_{red}$.
\end{defn}

\begin{defn}[Corollary 7.6 in \cite{marshall2006real}]\label{defn:mrrealreduced}
 A multiring $A$ is said to be \textit{real reduced} if the following properties holds for all 
$a,b,c,d\in F$:
\begin{enumerate}[i -]
 \item $1\ne0$;
 \item $a^3=a$;
 \item $c\in a+ab^2\Rightarrow c=a$;
 \item $c\in a^2+b^2$ and $d\in a^2+b^2$ implies $c=d$.
\end{enumerate}


A morphism of real reduced multirings is just a morphism of multirings. The category of real 
reduced multirings will be denoted by $\mathcal{MR}_{red}$.
\end{defn}

\section{Promiscuously Quadratic Ring}
 \hspace*{\parindent}
 Now, our interesting is analyze the Marshall's quotient in the case $(A,S)$ where $A$ is a commutative 
ring. Our subject is creat a new axiomatic for abstract quadratic forms, correlated to the theory of 
special groups \cite{dickmann2000special} and realsemigroups \cite{dickmann2004real}. Let start with an 
example:

\begin{ex}
 Let $A=\mathbb Z/_m(\mathbb Z^2\setminus\{0\})$ and $B=\mathbb Q/_m(\mathbb Q^2\setminus\{0\})$. Then
 $\overline a\in\overline b+\overline c$ (in $A$) if there exist $r,s,t\in\mathbb{Z}^*$ such that 
$ar^2=bs^2+ct^2$. As $r\ne0$, 
$$a=b\cdot\dfrac{s^2}{r^2}+c\cdot\dfrac{t^2}{r^2}$$
in $\mathbb Q$. Therefore $\overline a\in\overline b+\overline c$ in $A$ if and only if $\overline 
a\in\overline b+\overline c$ in $B$. Hence $B=M(G(\mathbb Q))$, we have that $A$ is a pre-special 
subgroup of $G(\mathbb Q)$.

Finally, $A=\{\overline p:p\mbox{ is free of squares in }\mathbb Z\}$.
\end{ex}

Motivated by this example, we developed this definition:

\begin{defn}[Promiscuously Quadratic Ring]\label{defn:pqring}
  A {\em promiscuously quadratic multiring (pq-multiring)} is a multiring $R$ satisfying the following 
properties:
\begin{description}
 \item [PQt -] (ternary) For all $x\in R$, $x^3=x$;
 \item [PQh -] (hyperbolic) For all $x,y\in R^\times$,  $y \in x + (-x)$.
\end{description}
\end{defn}
 
A pq-{\em multifield} $R$ will be called {\em pre-special} when holds
\begin{description}
 \item [PQfps-] (pre-special-mf) (i) If $ab = cd$ and $a \in c + d$ then $c \in a + b$ for all $a, b, c, 
d \in R^\times$; (ii) If $ab = cd = ef$, 
$a \in c+d$ and $c \in e+f$ then $a \in e+f$ for all $a, b, c, d, e, f \in R^\times$.
\end{description}

The notion of pq-multiring can be presented alternatively as a {\em first-order} structure in the 
language $L_{MR}$,  satisfying a (finite) set of $L_{MR}$-sentences of the form $\forall \vec{x} 
(\psi_1(\vec{x}) \rightarrow \psi_2(\vec{x}))$, where $\psi_i(\vec{x})$ is a pp-formula. The notion of 
real reduced pq-multiring  can be axiomatized by a set of {\em Horn geometric sentences}, i.e., 
sentences that are of the above form {\em or is a negation of an atomic formula}. The notion of 
pre-special multifield  can be axiomatized by a set of {\em Horn geometric sentences} and by the axiom 
$\forall x (x \neq 0 \rightarrow \exists y (x.y =1))$.

Considering the notions and results above and the development in \cite{ribeiro2016functorial}, the 
notions of (pre)special pq-multifield and  of (pre)special multifield coincide. The following result 
shows the coincidence between the above defined notion of real reduced pq-multiring and of real reduced 
multiring as presented in \cite{marshall2006real} (see Corollary 7.6).

\begin{prop} 
Let $M$ be a multiring satisfying [PQt] and [PQrr], then $M$ satisfies [PQh].
\end{prop}

The following result follows directly from the axiomatizations presented.

\begin{prop}
 The subcategory  $\PQM$  (respectively $\PQMrr$) is closed in  $\mathcal{MR} \subseteq L_{MR}-Str$  
under direct inductive limits and under arbitrary  products (respectively, non-empty products). Thus they 
are closed under reduced products (respectively, non trivial reduced products).
\end{prop}

\begin{defn}[Promiscuously Quadratic Pair]\label{defn:pqpair}
 A {\em promiscuously quadratic pair (pq-pair)} is a pair $(A,S)$ where $A$ is a ring, $S\subseteq 
A$, satisfying the following properties:
\begin{description}
 \item [pPQu -] (unity) $1\in S$; 
 \item [pPQm -] (multiplicative) $S\cdot S\subseteq S$;
 \item [pPQt -] (ternary) For all $a\in A$, there exist $r,s\in S$ such that $a^3r=as$;
 \item [pPQh -] (hyperbolic) For all $a,b\in A^\times$, there exist $r,s,t\in S$ such that $br=as-at$.
\end{description}

The notion of pq-pair can be presented alternatively as a {\em first-order} structure in language 
$L_{pMR}$, obtained from $L_{MR}$ by the addition of a unary predicate,  satisfying a (finite) set of 
$L_{pMR}$-sentences of the form $\forall \vec{x} (\psi_1(\vec{x}) \rightarrow \psi_2(\vec{x}))$, where 
$\psi_i(\vec{x})$ is a pp-formula.

A pq-pair $(A,S)$ will be called {\em real reduced} when  holds

[pPQrr -] (real reduced) (i) $0 \notin S$; (ii) $S+S\subseteq S$.

\end{defn}


\begin{prop}\label{pqred} 
If $(A,S)$ is pair  given by a  multiring and a proper pre-order in $A$, then 
$(A,S)$ is a  real reduced pq-pair.
\end{prop}

%
%


%

\begin{defn}[Morphism]\label{defn:pair-morph}
 A pq-pair morphism is just a $L_{pMR}$-homomorphism between pq-pairs. In more details, a map 
$f:(A,S)\rightarrow(B,T)$ will be a {\em morphism} of pq-pairs if $f:A\rightarrow B$ is an 
homomorphism of multirings and $f[S]\subseteq T$. 

The category of (real reduced) promiscuously quadratic pairs will be denoted by $\pPQM$ 
($\pPQMrr$).
\end{defn}

The following result follows directly from the axiomatizations presented.

\begin{prop}
 The subcategory  $\pPQM$  (respectively $\pPQMrr$) is closed in  $L_{pMR}-Str$  under direct inductive 
limits and under arbitrary  products (respectively, non-empty products). Thus they are closed under 
reduced products (respectively, non trivial reduced products).
\end{prop}

\begin{prop}
The mapping  $(A,S) \mapsto A/_mS$ is the object part of a  "definable" functor $q: \pPQM \rightarrow 
\PQM$. If $f:(A,S)\rightarrow(B,T)$ is a pq-pair morphism, then $\bar{f} : A/_mS \rightarrow B/_mT$, 
$\bar{f}([a]_S) = [f(a)]_T$, $a \in A$, is the action of $q$ on the arrow $f$.

The mapping $M \mapsto (M,\{1\})$ is the object part of a  "definable" functor $j: \PQM \rightarrow 
\pPQM$.

Moreover, $q \circ j \cong id_{\PQM}$.
\end{prop}

\begin{prop}
 The functors $q$ and $j$ preserve  products and direct inductive limits.
\end{prop}


\section{Quadratic form theory of pq-pairs}
 \hspace*{\parindent}
Let $(A,S)$ be a pq-pair.

(i)  A n-dimensional {\em form} over $(A,S)$ is just an n-tuple $\varphi=\langle\overline 
a_1,...,\overline a_n\rangle\in 
A$, $n \in \mathbb{N}$. The unique 0-dimensional form will be denoted by $\langle \ \rangle$.

(ii) - The {\em 0-ary isometry relation} over $(A,S)$ is the identity relation  $\langle \ \rangle = 
\langle \ \rangle$.\\
- The {\em 1-ary isometry relation} over $(A,S)$ is the equality relation in the associated pq-multiring 
$A/_mS$, i.e $\langle a\rangle\equiv_1\langle c\rangle$ iff $\bar{a} = \bar{c}$ in $A/_mS$.\\
- The {\em 2-isometry relation} over $(A,S)$ is the relation $\langle a,b\rangle\equiv_2\langle 
c,d\rangle$ if and 
only if the following equations holds in the associated pq-multiring $A/_mS$: $\bar{a}. \bar{b} = \bar{c} 
.\bar{d}$ and $\bar{a} +\bar{b} = \bar{c} +\bar{d}$. 
- For $n \geq 3$, the {\em n-isometry relation} is {\em transitive closure} of the relation $\langle 
a_1,..., a_n\rangle \approx \langle b_1,..., b_n\rangle \in A$ iff $\langle\overline a_1,...,\overline 
a_n\rangle = \langle\overline b_1,...,\overline b_n\rangle\in 
A/_mS$ {\bf or} $\exists 1\leq i, j, i', j' \leq n$, with $i \neq j$ and $i' \neq j'$, with $\langle 
a_i,a_j\rangle\equiv_2\langle b_{i'},b_{j'}\rangle$ and the corresponding $n-2$ forms obtained by erasing 
$a_i, a_j$ in the first member and $b_{i'}, b_{j'}$ in the second member, are $n-2$ isometric in $(A,S)$.

(Compare the present definition  with the definition 2.17 of T-isometry in \cite{dickmann2015faithfully} 
and with the classical Witt equivalence relation)

Note that the, for each $n \in \mathbb{N}$, the n-isometry relation is an equivalence relation in $A^n$. 

 If the extension of binary isometry to an $n$-ary relation, given by the standard 
construction in \cite{dickmann2000special}, is a transitive relation, then the  {n-isometry relation} 
over $(A,S)$ coincides with the with such the standard construction.




%
%
%
%
%
%
 
Considering the notions and results above and the development in \cite{ribeiro2016functorial}, the 
category of {\em reduced} real semigroups is equivalent to the category of {\em real reduced} 
pq-multirings. 

%

\section{Strictfication of pq-pairs}
Considering the notions and results above and the development in \cite{ribeiro2016functorial}, we 
formulate the following

{\bf Main question:} Is every pq-multiring $R$ isomorphic to $q(A,S)$ for some {\em univalent} pq-pair 
$(A,S)$ (i.e $A$ is a ring)?

If the answer for this question is positive, then there is, in particular, a concrete (and "definable") 
solution for the representation problem for special groups and real semigroups.

 \bibliographystyle{plain}
\bibliography{rings_many_squares}
\end{document}